\documentclass{amsart}
\usepackage{amsfonts}

\theoremstyle{plain}

\numberwithin{equation}{section}
\input{tcilatex}

\begin{document}
\title[On the coefficients of an expansion of $(1+1/x)^{x}$ related to
carleman's inequality]{On the coefficients of an expansion of $(1+1/x)^{x}$
related to carleman's inequality}
\author{Yue Hu}
\address{ College of Mathematics and Informatics, Henan Polytechnic
University, \\
Jaozuo City, Henan 454000, P.R. China.}
\email{huu3y2@163.com}
\author{Cristinel Mortici}
\address{ Valahia University of T\^{a}rgovi\c{s}te, Department of
Mathematics, Bd. Unirii 18, 130082 T\^{a}rgovi\c{s}te, Romania }
\email{cmortici@valahia.ro}

\begin{abstract}
In this note, we present new properties for the sequence $b_1=\frac{1}{2},$ $%
b_n=\frac{1}{n}(\frac{1}{n+1}-\sum_{k=0}^{n-2}\frac{b_{n-k-1}}{k+2}),n\geq2,$
arising in some refinements of Carleman's inequality. Our results extend
some results of Yang [Approximations for constant $e$ and their applications
J. Math.  Anal. Appl. 262 (2001) 651-659] and Alzer and Berg [some classes
of completely  monotonic functions Ann. Acad. Sci. Fennicae 27(2002)
445-460].
\end{abstract}

\maketitle



\vspace{.08in} \noindent \textbf{Keywords}: constant $e$; sequence;
inequalities; integral representation 

\section{Introduction}

The following Carleman inequality \cite{carl}%
\begin{equation*}
\sum_{n=1}^{\infty }\left( a_{1}a_{2}\cdots a_{n}\right)
^{1/n}<e\sum_{n=1}^{\infty }a_{n},
\end{equation*}%
whenever $a_{n}\geq 0,$ $n=1,2,3,\ldots ,$ with $0<\sum_{n=1}^{\infty
}a_{n}<\infty ,$ has attracted the attention of many authors in the recent
past. We refer for example to \cite{hu}, or to the work of Yang \cite{Yang},
who proved 
\begin{equation*}
\sum_{n=1}^{\infty }\left( a_{1}a_{2}\cdots a_{n}\right)
^{1/n}<e\sum_{n=1}^{\infty }\left( 1-\sum_{k=1}^{6}\frac{b_{k}}{\left(
n+1\right) ^{k}}\right) a_{n}
\end{equation*}%
with $b_{1}=1/2,$ $b_{2}=1/24,$ $b_{3}=1/48,$ $b_{4}=73/5760$, $%
b_{5}=11/1280,$ $b_{6}=1945/580608.$ In the final part of his paper, Yang 
\cite{Yang} conjectured that if 
\begin{equation*}
\left( 1+\frac{1}{x}\right) ^{x}=e\left( 1-\sum_{n=1}^{\infty }\frac{b_{n}}{%
\left( x+1\right) ^{n}}\right) ,\ \ \ x>0,  \label{S0}
\end{equation*}%
then $b_{n}>0,$ $n=1,2,3,\ldots $ .\newline
\indent Later , this conjecture was proved and discussed by Yang \cite{yang}%
, Gylletberg and Yan \cite{Yan}, and Yue \cite{Yue}, using the recurrence 
\begin{equation}
b_1=\frac{1}{2}, b_n=\frac{1}{n}(\frac{1}{n+1}-\sum_{k=0}^{n-2}\frac{%
b_{n-k-1}}{k+2})\qquad(n=2,3,\dotsm).
\end{equation}
Here we observe that it is easy to calculate the value of $b_n$ by computer
programs, but it is hard to extract any properties about $b_n$ from (1.1).
Also the proofs of $b_n>0$ provided were quite complicated.\newline
\indent In this note, we use an argument of Alzer and Berg \cite{Alzer} to
derive an integral representation of $b_n$ and then we get some new
properties about $b_n$ at once. 

\section{Lemmas}

In order to prove our main results we need the following two lemmas. A proof
of Lemma 1 is given in \cite{Alzer} .\newline
\textbf{Lemma 1.}\emph{\qquad }Let 
\begin{equation}
f(x)=(x+1)[e-(\frac{1+x}{x})^{x}]\qquad ( x>0 ).
\end{equation}
Then we have 
\begin{equation}
f(x)=\frac{e}{2}+\frac{1}{\pi}\int_{0}^{1}\frac{s^s(1-s)^{1-s}sin(\pi{s})}{%
x+s}ds .
\end{equation}
The next lemma can be found in [5].\newline
\textbf{Lemma 2.}\emph{\qquad }Let $h(s) $ is a continuous function on
[0,1]. Then 
\begin{equation}
\lim_{n \rightarrow \infty}n\int_{0}^{1}{s^n}h(s)ds=h(1) .
\end{equation}

\section{Results}

The main result of this paper which incorporate the announced new properties
of the sequence ${\{b_n\}_{n\geq1}}$ in the following. \newline
\textbf{Theorem 1.}\emph{\qquad }Let $\{b_n\}_{n\geq1}$be the sequence
defined by (1.1) , and let 
$$
g(s)=\frac{1}{\pi}{s^s}(1-s)^{1-s}sin(\pi{s}).%
$$
Then 
\begin{equation}
b_n=\frac{1}{e}\int_{0}^{1}g(s)(1-s)^{n-2}ds=\frac{1}{e}%
\int_{0}^{1}g(s)s^{n-2}ds,
\end{equation}
\begin{equation}
0<b_n\leq\frac{1}{n(n+1)}\qquad(n=1,2,\dotsm),
\end{equation}
\begin{equation}
b_{n+1}<b_{n} \qquad(n=1,2,\dotsm),
\end{equation}
\begin{equation}
\lim_{n\rightarrow\infty}\frac{b_{n+1}}{b_{n}}=1.
\end{equation}

\begin{proof}
Let 
\begin{equation*}
t=\frac{1}{1+x}\qquad(x>0),
\end{equation*}
\begin{equation}
f(t)=e-(1-t)^{1-{1/t}}.
\end{equation}
From (1.1) we obtain 
\begin{equation}
b_{n}=\frac{f^{(n)}(0)}{n!e}\qquad(n=1,2,\dotsm).
\end{equation}
On the other hand , by Lemma 1 , we have 
\begin{equation}
f(t)=\frac{e}{2}t+\int_{0}^{1}g(s)\frac{t^2}{1+(s-1)t}ds ,
\end{equation}
where ${0}\leq{t}\leq{1}$ .\newline
Differentiation gives 
\begin{equation*}
f^{\prime}(0)=\frac{e}{2} ,
\end{equation*}
\begin{equation}
f^{(n)}(0)=n!\int_{0}^{1}g(s)(1-s)^{n-2}ds\qquad (n=2,3,\dotsm) .
\end{equation}
So that (3.6) and (3.8) yield 
\begin{equation*}
b_1=\frac{1}{2} ,
\end{equation*}
\begin{equation}
b_{n}=\frac{1}{e}\int_{0}^{1}g(s)(1-s)^{n-2}ds =\frac{1}{e}%
\int_{0}^{1}g(s)s^{n-2}ds \qquad(n=2,3,\dotsm),
\end{equation}
where the last step uses the change of variable $s=1-t.$ Which prove (3.1).%
\newline
The relation (3.2) and (3.3) follow immediately from (3.9) and (1.1) .%
\newline
Now , we prove (3.4). Since $g(1)=0$ , we cannot directly use (3.9) to
calculating this limit by Lemma 2.\newline
Integrating by parts , we find from (3.9) 
\begin{equation}
b_n=-\frac{1}{(n-1)e}\int_{0}^{1}h(s)s^{n-1}ds\qquad(n=2,3,\dotsm),
\end{equation}
where 
\begin{equation*}
h(s)={s^s}(1-s)^{1-s}[cos(\pi{s})-\frac{sin(\pi{s})}{\pi}\ln\frac{1-s}{s}].
\end{equation*}
Noting 
\begin{equation*}
\lim_{s\rightarrow0^{+}}\frac{sin(\pi{s})}{\pi}\ln\frac{1-s}{s}=0,
\end{equation*}
\begin{equation*}
\lim_{s\rightarrow1^{-}}\frac{sin(\pi{s})}{\pi}\ln\frac{1-s}{s}=0,
\end{equation*}
from Lemma 2 and (3.9) , we have%
\begin{equation*}
\begin{split}
&\quad\lim_{n\rightarrow\infty}\frac{b_{n+1}}{b_{n}} \\
&=\lim_{n\rightarrow\infty}\frac{(n-1)\int_{0}^{1}h(s)s^{n}ds}{%
n\int_{0}^{1}h(s)s^{n-1}ds} \\
&=\lim_{n\rightarrow\infty}\frac{(n-1)^{2}}{n^{2}}\frac{(n)%
\int_{0}^{1}h(s)s^{n}ds}{(n-1)\int_{0}^{1}h(s)s^{n-1}ds} \\
&=\frac{h(1)}{h(1)}=\frac{-1}{-1}=1.
\end{split}%
\end{equation*}%
This completes the proof of Theorem 1.
\end{proof}

\indent Remark That relation (3.1) can be rewritten as follows: 
\begin{equation*}
\int_{0}^{1}g(s)(s)^{n-2}ds=\int_{0}^{1}g(s)(1-s)^{n-2}ds=eb_{n}%
\qquad(n=2,3,\dotsm).
\end{equation*}%
\newline
which means that they are extension of [1,Lemma2]. For example 
\begin{equation*}
\int_{0}^{1}g(s)ds=eb_{2}=\frac{e}{24},\qquad \int_{0}^{1}g(s)sds=eb_{3}=%
\frac{e}{48},
\end{equation*}
and 
\begin{equation*}
\begin{split}
&\quad \int_{0}^{1}\frac{1}{s}g(s)ds=\int_{0}^{1}\frac{1}{1-s}g(s)ds \\
&=\int_{0}^{1}(1+s+s^{2}+\dotsm)g(s)ds \\
&=e\sum_{n=2}^{\infty}b_{n} =e\sum_{n=1}^{\infty}b_{n}-eb_{1} \\
&=e(1-\frac{1}{e})-\frac{e}{2} =\frac{e}{2}-1.
\end{split}%
\end{equation*}

\end{document}